\font\we=cmb10 at 14.4truept
\font\li=cmb10 at 12truept
\noindent
\centerline {\we Riemann-Roch, Stability and}
\centerline {\we New Non-Abelian Zeta
Functions for Number Fields}
\vskip 0.45cm
\centerline {\li Lin WENG}
\vskip 1.0cm
\noindent
In this paper, we introduce a geometrically stylized arithmetic cohomology for
number fields. Based on such a cohomology, we define and study
new yet genuine non-abelian zeta functions for number fields, using an
intersection stability.
\vskip 0.30cm
\centerline {\li 1. A New Cohomology}
\vskip 0.30cm
\noindent
{\bf (1.1) Arithmetic Cohomology Groups}
\vskip 0.30cm
Let $F$ be a number field with discriminant $\Delta_F$. Denote its
(normalized) absolute values by
$S_F$, and write $S_F=S_{\rm\,  fin}\buildrel \cdot\over\cup S_\infty$,
where $S_\infty$ denotes the collection of all archimedean valuations.
For simplicity, we use $v$ (resp. $\sigma$) to denote elements in
$S_{\rm\,  fin}$ (resp. $S_\infty$).

Denote by ${\bf A}={\bf A}_F$ the ring of adeles of $F$, by ${\rm\,
Gl}_r({\bf A})$ the rank $r$ general linear group over {\bf A}, and write
${\bf A}:={\bf A}_{\rm\,  fin}\oplus {\bf A}_\infty$ and ${\rm\,  GL}_r({\bf
A}):={\rm\,  GL}_r({\bf A})_{\rm\,  fin}\times {\rm\,  GL}_r({\bf A})_\infty$
according to their finite and infinite parts.

For any $g=(g_{\rm\,  fin}:g_\infty)=(g_v;g_\sigma)\in {\rm\,  GL}_r({\bf A})$,
define the injective morphism $i(g):=i(g_\infty):F^r\to
{\bf A}^r$ by $(f)\mapsto (f;g_\sigma\cdot f)$. Let $F^r(g):=
{\rm\,  Im}\big(i(g)\big)$ and set
$${\bf A}^r(g):=\{(a_v;a_\sigma)\in {\bf A}^r:g_v(a_v)\in {\cal
O}_v^r,\ \forall v;\ {\rm\,  and}\
\exists f\in F^r\ {\rm\,  s.t.}\ g_v(f)\in {\cal
O}_v^r,\,\forall v\ {\rm\,  and}\ (f;a_\sigma)=i(g_\infty)(f)\}.$$ Thenwe
have the following 9-diagram with exact columns and rows:
$$\matrix{&&0&&0&&0&&\cr
&&\downarrow&&\downarrow&&\downarrow&&\cr
0&\to&{\bf A}^r(g)\cap F^r(g)&\to&{\bf A}^r(g)&\to&{\bf A}^r(g)/
{\bf A}^r(g)\cap F^r(g)&\to&0\cr
&&\downarrow&&\downarrow&&\downarrow&&\cr
0&\to& F^r(g)&\to&{\bf A}^r&\to&{\bf A}^r/F^r(g)&\to&0\cr
&&\downarrow&&\downarrow&&\downarrow&&\cr
0&\to& F^r(g)/{\bf A}^r(g)\cap F^r(g)&\to&{\bf A}^r/{\bf A}^r(g)&\to&
{\bf A}^r/{\bf A}^r(g)+F^r(g)&\to&0\cr
&&\downarrow&&\downarrow&&\downarrow&&\cr
&&0&&0&&0.&&\cr}$$

Motivated by this and Weil's adelic cohomology theory for divisors over
algebraic curves, (see e.g., [24] and [29]), we introduce the
following
\vskip 0.30cm
\noindent
{\bf Definition.} {\it For any $g\in {\rm\,  GL}_r({\bf A})$, define its
0-th and 1-st arithmetic cohomology groups by
$$H^0({\bf A}_F,g):={\bf A}^r(g)\cap F^r(g),\qquad{\rm\,  and}\qquad 
H^1({\bf A}_F,g):=
{\bf A}^r/{\bf A}^r(g)+F^r(g).$$}

\noindent
{\bf Theorem.} (Serre Duality=Pontrjagin Duality) {\it As locally compact
groups, $$H^1({\bf A}_F,g)\simeq \widehat{H^0({\bf A}_F,k_F\otimes 
g^{-1})}.$$ Here $k_F$
denotes an idelic dualizing element of $F$, and $\hat\cdot$  the
Pontrjagin dual. In particular,
$H^0$ is  discrete and $H^1$ is compact.}
\vskip 0.30cm
\noindent
{\it Remark.} For $v\in S_{\rm\,  fin}$, denote by $\partial_v$ the local
different of $F_v$, the $v$-completion of $F$ at $v$, and by
${\cal O}_v$ the valuation ring with $\pi_v$ a local parameter. Then
$\partial_v=:\pi_v^{{\rm\,  ord}_v(\partial_v)}\cdot {\cal O}_v$. We call
$\kappa_F:=(\partial_v^{{\rm\,  ord}_v(\partial_v)};1)\in {\bf I}_F:={\rm\, 
GL}_1({\bf A})$
an idelic dualizing element of $F$.
\vskip 0.30cm
\noindent
{\it Proof.} As usual, introduce a basic character $\chi$ on {\bf A} by
$\chi=:(\chi^{(r)}_v;\chi^{(r)}_\sigma)$ where $\chi_v:=\lambda_v\circ
{\rm\,  Tr}^{F_v}_{{\bf Q}_v}$ with $\lambda_v:{\bf Q}_v\to {\bf Q}_v/
{\bf Z}_v\hookrightarrow {\bf Q}/{\bf Z}\hookrightarrow {\bf R}/{\bf Z}$, and
$\chi_\sigma:=\lambda_\infty\circ {\rm\,  Tr}_{\bf R}^{F_\sigma}$ with
$\lambda_\infty:{\bf R}\to {\bf R}/{\bf Z}$.  Then the pairing
$(x,y)\mapsto e^{2\pi i\chi(x\cdot y)}$ induces natural isomorphisms
$\widehat {{\bf A}^r}\simeq {\bf A}^r$ (as
locally compact groups) and $(F^r)^\perp\simeq F^r$ (as discrete subgroups).
With this, a direct local calculation shows that
  $\big({\bf A}^r(g)\big)^\perp={\bf A}^r(\kappa_F\otimes g^{-1})$ and
$\big(F^r(g)\big)^\perp\simeq F^r(\kappa_F\otimes g^{-1})$. This
completes the proof since
$$\Big({\bf A}^r(g)\cap F^r(g)\Big)^\perp=\big({\bf A}^r(g)\big)^\perp
+\big(F^r(g)\big)^\perp.$$
\vskip 0.30cm
\noindent
{\bf (1.2) Arithmetic Counts}
\vskip 0.30cm
Motivated by the Pontrjagin duality and the fact that the dimension of a 
vector space is
equal to the dimension of its dual, one  basic principal we
adopt in counting locally compact groups is the following:
\vskip 0.30cm
\noindent
{\bf Counting Axiom.} {\it If $\#_{\rm\,  ga}$  counts a certain
class of locally compact groups $G$,  then $\#_{\rm\,  ga}(G)=
\#_{\rm\,  ga}(\hat G)$.}
\vskip 0.30cm
Practically, our counts of arithmetic cohomology groups are based on the 
Fourier inverse
formula, or more accurately, the Plancherel formula in Fourier analysis 
overlocally compact
groups. (See e.g. [9].)

While  any reasonable test function on ${\bf A}^r$ would do,  as a
continuation of a more classical mathematics and also for simplicity, we 
set$f:=\prod_v
f_v\cdot\prod_\sigma f_\sigma$. Here $f_v$ is the characteristic function 
of${\cal O}_v^r$;
$f_\sigma(x_\sigma):=e^{-\pi|x_\sigma|^2/2}$ if $\sigma$ is real; and
$f_\sigma(x_\sigma):=e^{-\pi|x_\sigma|^2}$ if $\sigma$ is complex.
Moreover, we take the following  normalization for the Haar measure
$dx$, which we call {\it standard}, on {\bf A}: locally for $v$, $dx$ is
the measure for which
${\cal O}_v$ gets the measure $N(\partial_v)^{-1/2}$, while for
$\sigma$ real (resp. complex), $dx$ is the ordinary Lebesgue measure
(resp. twice the ordinary Lebesgue measure).
\vskip 0.30cm
\noindent
{\bf Definition.} {\it (1) The arithmetic counts of the 0-th and the
1-st arithmetic cohomology groups for  $g\in {\rm\,  GL}_r({\bf A})$ are
defined to be
$$\eqalign{\#_{\rm\,  ga}\big(H^0({\bf A}_F,g)\big):=&\#_{\rm\, 
ga}\Big(H^0({\bf A}_F,g); f,
dx\Big):=\int_{H^0({\bf A}_F,g)}|f(x)|^2dx;\cr
\#_{\rm\,  ga}\big(H^1({\bf A}_F,g)\big):=&\#_{\rm\,  ga}\Big(H^1({\bf 
A}_F,g); \hat f,
d\xi\Big):=\int_{H^1({\bf A}_F,g)}|\hat f(\xi)|^2d\xi.\cr}$$ Here $dx$ denotes
(the restriction of) the standard Haar measure on {\bf A},  $d\xi$
(the induced quotient measure from) the dual measure (with
respect to $\chi$), and
$\hat f$  the corresponding Fourier transform of $f$;
\vskip 0.30cm
\noindent
(2) The 0-th and the 1-st arithmetic cohomologies of
$g\in {\rm\,  GL}_r({\bf A})$ are defined to be
$$h^0({\bf A}_F,g):=\log \Big(\#_{\rm\,  ga}\big(H^0({\bf 
A}_F,g)\big)\Big)\qquad
{\rm\,  and}\qquad h^1({\bf A}_F,g):=\log \Big(\#_{\rm\,  ga}\big(H^1({\bf 
A}_F,g)\big)\Big).$$}

\noindent
{\bf (1.3) Serre Duality and Riemann-Roch}
\vskip 0.30cm
For the arithmetic cohomologies just introduced, we have the following
\vskip 0.30cm
\noindent
{\bf Theorem.} (1) (Serre Duality) $h^1({\bf A}_F,g)=
h^0({\bf A}_F,\kappa_F\otimes g^{-1})$;
\vskip 0.30cm
\noindent
(2) (Riemann-Roch Theorem)
$$h^0({\bf A}_F,g)-h^1({\bf A}_F,g)={\rm\,  deg}(g)-{r\over 
2}\cdot\log|\Delta_F|.$$

\noindent
{\it Proof.} By the choice of our $f$ in the definition, (1) is a direct 
consequence of the
topological  Serre duality, i.e., Theorem 1.1, and the Plancherel Formula, 
while
(2) is a direct consequence of the Serre duality just proved and Tate's
Riemann-Roch theorem ([28, Thm. 4.2.1] and/or [14, XIV, \S6]), i.e, the 
Poisson summation formula,
by the fact that
$\Big(H^0({\bf A}_F,g)\Big)^\perp =H^0({\bf A}_F,\kappa_F\otimes g^{-1})$.
This then completes the proof.

Often, for our own convenience, we also write $e^{h^i({\bf A}_F,g)}=
\#_{\rm\,  ga}\big(H^i({\bf A}_F,g)\big)$ simply as $H^i_{\rm\,  ga}(F,g),\ 
i=1,2$. With this, the above additive
version may be rewritten as
\vskip 0.30cm
\noindent
{\bf Theorem}$'$. (1) (Serre Duality) $H^1_{\rm\,  ga}(F,g)=H^0_{\rm\, ga}
(F,\kappa_F\otimes g^{-1})$;
\vskip 0.30cm
\noindent
(2) (Riemann-Roch Theorem)
$H^0_{\rm\,  ga}(F,g)=H^1_{\rm\,  ga}(F,g)\cdot N(g)\cdot 
N(\kappa_F)^{-{r\over 2}}$, where as
usual $N(g)$ denotes $e^{{\rm\,  deg}(g)}$.
\vskip 0.30cm
\noindent
{\it Remarks.} (1) Our work here is motivated by the works of Weil, Tate,
van der Geer-Schoof, and Li, as well as the works of Lang, Arakelov,
Szpiro, Parshin,  Moreno, Neukirch,  Deninger, Connes, and  Borisov. For 
details,
please see the references below, in particular [31]. Also, it would be 
extremely
interesting if one could relate the work here with that of Connes [5] and 
Deninger [7, 8].

\noindent
(2)  One may apply the discussion in this paper to  wider classes of 
(multiplicative) characters
and test functions. We leave this to the reader. (See e.g., [28], [29] and 
[18].)
\vskip 0.8cm
\centerline {\li 2. New Non-Abelian Zeta Functions}
\vskip 0.30cm
\noindent
{\bf (2.1) Intersection Stability}
\vskip 0.30cm
For a metrized vector sheaf $({\cal E},\rho)$ on ${\rm\,  Spec}({\cal 
O}_F)$, define its
associated
$\mu$-invariant by $$\mu({\cal E},\rho):={{{\rm\,  deg}_{\rm\,  Ar}({\cal 
E},\rho)}\over
{{\rm\,  rank}({\cal E})}},$$ where ${\cal O}_F$ denotes the ring of 
integers of a number field $F$ and
${\rm\,  deg}_{\rm\,  ar}$  the Arakelov degree of $({\cal E},\rho)$. (See 
e.g. [16].) By
definition, a proper sub metrized vector sheaf $({\cal E}_1,\rho_1)$ of 
$({\cal E},\rho)$
consists of a proper sub vector sheaf ${\cal E}_1$ of ${\cal E}$ such that 
$\rho_1$
is induced from the restriction of $\rho$ via the injection ${\cal 
E}_1\hookrightarrow {\cal E}$.
\vskip 0.30cm
\noindent
{\bf Definition.} {\it A metrized vector sheaf $({\cal E},\rho)$ is called 
stable (resp.
semi-stable) if for all proper sub metrized vector sheaf $({\cal 
E}_1,\rho_1)$ of
$({\cal E},\rho)$,
$$\mu({\cal E}_1,\rho_1)\,<\,\mu({\cal E},\rho)\qquad(resp.\qquad
\mu({\cal E}_1,\rho_1)\leq\mu({\cal E},\rho)).$$}

\noindent
{\it Remarks.} (1) Despite the fact that we define it independently,
the intersection stability in arithmetic, motivated by Mumford's work [20] 
in geometry, was first
introduced by Stuhler in [26,27],  see also [11,12,19 and 4]. Standard
facts concerning Harder-Narasimhan filtrations and Jordan-H\"older graded 
metrized vector sheaves
hold in this setting as well. For details, see e.g., [4,19 and 31].

\noindent
(2) The intersection stability plays a key role in our work on non-abelian 
class field theory
for Riemann surfaces in [30]. Motivated by this, as a fundamental problem, 
we ask
whether a Narasimhan-Seshadri type correspondence holds in arithmetic in [31].
\vskip 0.30cm
For $g=(g_v;g_\sigma)\in {\rm\,  GL}_r({\bf A})$, introduce a torsion-free 
${\cal O}_F$-module
$$H^0({\bf A}_F,g)_{\rm\,  fin}:=H^0({\rm\,  Spec}({\cal O}_F),g):=\{f\in 
F^m:g_vf\in
{\cal O}_v^r,\forall v\}$$ in
$F^r$. Denote the associated vector sheaf on ${\rm\,  Spec}({\cal O}_F)$ by 
${\cal E}(g)$, that is,
$${\cal E}(g):=\widetilde {H^0({\rm\,  Spec}({\cal O}_F),g)}.$$
Moreover, note that $F^r$, via completion, is densely embedded in ${\bf 
A}_\infty^r$. Thus to
introduce metrics on ${\cal E}(g)$ is the same as to assign metrics on the 
determinants, i.e., on
the top exterior products, of the associated data. (See e.g., [16, Chap. 
V].) Hence without loss of
generality, we may assume that $r=1$. In this case, the metric on ${\cal 
E}(g)$ associated to $g$
is defined to be the one such that for the rational section $1\in F$,
$$\|1\|_\sigma:=\|g_\sigma\|_\sigma:=|g_\sigma|^{N_\sigma:=[F_\sigma:{\bf 
Q}_\sigma]}.$$ Denote
such a metric on ${\cal E}(g)$ by $\rho(g)$ for $g\in {\rm\,  GL}_r({\bf A})$.

As such, we obtain a canonical map $({\cal E}(\cdot),\rho(\cdot)):{\rm\, 
GL}_r({\bf A})\to
\Omega_{{\rm\,  Spec}({\cal O}_F),r}$ by assigning $g$ to $({\cal 
E}(g),\rho(g))$, where
$\Omega_{{\rm\,  Spec}({\cal O}_F),r}$ denotes the collection of all 
metrized vector sheaves of
rank $r$ over ${\rm\,  Spec}({\cal O}_F)$. Clearly, $({\cal 
E}(\cdot),\rho(\cdot))$ factors through
the quotient group ${\rm\,  GL}_r(F)\backslash {\rm\,  GL}_r({\bf A})$ 
where ${\rm\,  GL}_r(F)$ is
embedded diagonally in ${\rm\,  GL}_r({\bf A})$. Denote this resulting map  by
$({\cal E}(\cdot),\rho(\cdot))$ too by an abuse of notation.

Denote by ${\cal M}_{F,r}(d)$ the subset of $\Omega_{{\rm\,  Spec}({\cal 
O}_F),r}$ consisting
of semi-stable metrized vector sheaves of (Arakelov) degree $d$. Since for 
afixed degree,
the semi-stability condition is a bounded and closed one, with respect to 
the natural topology,
${\cal M}_{F,r}(d)$ is compact. (See e.g. [11,26 and 27].)

Denote by ${\cal M}_{{\bf A}_F,r}(d)\subset {\rm\,  GL}_r(F)\backslash 
{\rm\,  GL}_r({\bf A})$ the
inverse image of ${\cal M}_{F,r}(d)$ with respect to $({\cal 
E}(\cdot),\rho(\cdot))$, and denote
the corresponding map by
$$\Pi_{F,r}(d):{\cal M}_{{\bf A}_F,r}(d)\to {\cal M}_{F,r}(d)$$ which we 
call the (algebraic)
moment map. As a subquotient of ${\rm\,  GL}_r({\bf A})$, ${\cal M}_{{\bf 
A}_F,r}(d)$ admits a
natural topology, the induced one. Moreover, by a general result due to 
Borel [2], which in
our case is more or less obvious, the fibers of $\Pi_{F,r}(d)$ are all 
compact.  Thus in
particular, ${\cal M}_{{\bf A}_F,r}(d)$, which we call the moduli space of 
semi-stable
adelic bundles of rank $r$ and degree $d$, is compact. In particular, as a 
subquotient of
${\rm\,  GL}_r({\bf A})$, ${\cal M}_{{\bf A}_{F,r}}(d)$ carries a natural 
measure induced from the
standard one on ${\rm\,  GL}_r({\bf A})$, which we call the {\it Tamagawa 
measure}, and denote  by
$d\mu_{{\bf A}_F,r}(d)$. For the same reason, there is also a natural 
measure on
${\cal M}_{F,r}(d)$, which we call the {\it hyperbolic measure}, and denote 
it by $d\mu_{F,r}(d)$.

Clearly, the total volumes of ${\cal M}_{{\bf A}_F,r}(d)$ (resp. ${\cal 
M}_{F,r}(d)$) with
respect to $d\mu_{{\bf A}_F,r}(d)$ (resp. $d\mu_{F,r}(d)$) are very 
important non-commutative
invariants for number fields.
Note that according to what we call the Bombieri-Vaaler trick [1], i.e., by 
multiplying $g$ with
$(1;e^{t_\sigma})$ where $t_\sigma:=N_\sigma\cdot t$ with
$t\in {\bf R}$, we obtain a natural isomorphism between ${\cal M}_{{\bf 
A}_F,r}(d)$ and
${\cal M}_{{\bf A}_F,r}(d-n\cdot t)$ where $n:=[F:{\bf Q}]$. (Even though 
itis an open problem that
semi-stability is closed under tensor operation [4], the case here in which 
one is of rank 1 is
rather obvious.) Consequently, the above volumes
are independent of degrees $d$. Denote them by $W_F(r)$ and $w_F(r)$
respectively.
\vskip 0.30cm
\noindent
{\bf (2.2) Functional Equation: A Formal Calculation}
\vskip 0.30cm
Let $F$ be a number field with discriminant $\Delta_F$. Denote by ${\cal 
M}_{{\bf A}_F,r}$ the
moduli space of semi-stable adelic bundles of rank $r$, that is, ${\cal 
M}_{{\bf
A}_F,r}:=\buildrel\cdot\over\cup_{N\in {\bf R}_+} {\cal M}_{{\bf A}_F,r}[N]$
where ${\cal M}_{{\bf A}_F,r}[N]:={\cal M}_{{\bf A}_F,r}(\log\,N)$.
By using the Bombieri-Vaaler trick in (2.1), as topological spaces,
${\cal M}_{{\bf A}_F,r}\simeq {\cal M}_{{\bf A}_F,r}[|\Delta_F|^{r\over 
2}]\times {\bf R}_+$.
Hence we obtain a natural measure $d\mu$ on ${\cal M}_{{\bf A}_F,r}$ from 
the Tamagawa
measures on ${\cal M}_{{\bf A}_F,r}[N]$ and ${{dt}\over t}$ on ${\bf R}_+$.

For any $E\in {\cal M}_{{\bf A}_F,r}$, define $H_{\rm\, 
ga}^i(F,E):=H_{\rm\,  ga}^i(F,g)
=e^{h^i({\bf A}_F,g)}$ for any $g\in {\rm\,  GL}_r({\bf A})$ such that 
$E=[g]$. Since for any
$a\in {\rm\,  GL}_r(F)$, $H_{\rm\,  ga}^i(F,a\cdot g)=H^i_{\rm\, ga}(F,g)$, 
$H_{\rm\,  ga}^i(F,E)$ is
well-defined for $i=0,\,1$.
\vskip 0.30cm
With respect to fixed real constants $A,B,C,\alpha$ and $\beta$, introduce 
the formal
integration $Z_{F,r;A,B,C;\alpha,\beta}(s)$ as follows:
$$Z_{F,r;A,B,C;\alpha,\beta}(s):=\big(|\Delta_F|^{-{{rB}\over 2}}\big)^s
\int_{E\in {\cal M}_{{\bf A}_F,r}}\Big(H_{\rm\,  ga}^0(F,E)^A\cdot
N(E)^{Bs+C}-N(E)^{\alpha s+\beta}\Big)d\mu(E).$$ Then formally,
$$Z_{F,r;A,B,C;\alpha,\beta}(s)=I(s)-II(s)+III(s),$$ where
$$\eqalign{I(s):=&\big(|\Delta_F|^{-{{rB}\over 2}}\big)^s
\int_{E\in {\cal M}_{{\bf A}_F,r},N(E)\leq |\Delta_F|^{r\over 
2}}\Big(H_{\rm\,  ga}^0(F,E)^A\cdot
N(E)^{Bs+C}-N(E)^{\alpha s+\beta}\Big)d\mu(E);\cr
II(s):=&\big(|\Delta_F|^{-{{rB}\over 2}}\big)^s
\int_{E\in {\cal M}_{{\bf A}_F,r},N(E)\geq |\Delta_F|^{r\over 
2}}N(E)^{\alpha s+\beta}d\mu(E);\cr
III(s):=&\big(|\Delta_F|^{-{{rB}\over 2}}\big)^s
\int_{E\in {\cal M}_{{\bf A}_F,r},N(E)\geq |\Delta_F|^{r\over 2}}H_{\rm\, 
ga}^0(F,E)^A\cdot
N(E)^{Bs+C} d\mu(E).\cr}$$
By Theorem$'$ 1.3, i.e., the multiplicative Serre duality and Riemann-Roch 
theorem, we
have
$$III(s)=\big(|\Delta_F|^{-{{rB}\over 2}}\big)^{-s-{{A+2C}\over B}}
\int_{E\in {\cal M}_{{\bf A}_F,r},N(E)\leq |\Delta_F|^{r\over 2}}H_{\rm\, 
ga}^0(F,E)^A\cdot
N(E)^{B(-s-{{A+2C}\over B})+C} d\mu(E),$$ by the fact that $N(E_1\otimes 
E_2^\vee)=
N(E_1)^{{\rm\,  rank}(E_2)}\cdot N(E_2)^{-{\rm\,  rank}(E_1)}.$ Hence, 
formally,
$$Z_{F,r;A,B,C;\alpha,\beta}(s)=I(s)+I(-s-{{A+2C}\over B})-II(s)+IV(s),$$ where
$$IV(s):=\big(|\Delta_F|^{-{{rB}\over 2}}\big)^{-s-{{A+2C}\over B}}
\int_{E\in {\cal M}_{{\bf A}_F,r},N(E)\leq |\Delta_F|^{r\over 2}}
N(E)^{\alpha\cdot(-s-{{A+2C}\over B})+\beta}d\mu(E).$$

Moreover, by definition,
$$\eqalign{-II(s)=&-\int_{E\in {\cal M}_{{\bf A}_F,r},
N(E)\geq
\Delta_F^{r\over 2}}
N(E)^{\alpha s+\beta}d\mu(E)
=-\int_{{\cal
M}_{{\bf A}_F,r}[|\Delta_F|^{r\over 2}]} d\mu(E)\cdot\int_1^\infty
t^{\alpha t+\beta}{{dt}\over t}\cr
=&-W_F(r)\cdot
{{t^{\alpha s+\beta}}\over
{\alpha s+\beta}}\Big|_1^\infty=W_F(r)\cdot {1\over {\alpha 
s+\beta}},\cr}$$provided that
$\alpha s+\beta<0$.

Similarly,
$$IV(s)=\int_{{\cal M}_{{\bf A}_F,r}[|\Delta_F|^{r\over 
2}]}d\mu(E)\cdot\int_0^1
t^{\alpha(-s-{{A+2C}\over B})
+\beta}{{dt}\over t}
=W_F(r)\cdot
{1\over {\alpha(-s-{{A+2C}\over B})
+\beta}}$$ provided that $\alpha(-s-{{A+2C}\over
B}) +\beta>0$.

Therefore, formally,
$$Z_{F,r;A,B,C;\alpha,\beta}(s)
=I(s)+I\big(-s-{{A+2C}\over B}\big)
+W_F(r)\cdot\Big({1\over
{\alpha s+\beta}}+{1\over {\alpha(-s-{{A+2C}\over
B})+\beta}}\Big).$$
As a direct consequence, we have the following
\vskip 0.30cm
\noindent
{\bf Functional Equation.} {\it With the same notation as above, formally,
$$Z_{F,r;A,B,C;\alpha,\beta}(s)=Z_{F,r;A,B,C;\alpha,\beta}(-s-{{A+2C}\over
B}).$$}

\noindent
{\bf (2.3) Non-Abelian Zeta Functions for Number Fields}
\vskip 0.45cm
To justify the arguments in (2.2), we consider convergences of two
types.
\vskip 0.30cm
\noindent
{\it Type 1}. Convergence for $II(s)$ and $IV(s)$, where
$$\eqalign{II(s)=&\big(|\Delta_F|^{-{{rB}\over 2}}\big)^s
\int_{E\in {\cal M}_{{\bf A}_F,r},N(E)\geq |\Delta_F|^{r\over 
2}}N(E)^{\alpha s+\beta}d\mu(E);\cr
IV(s)=&\big(|\Delta_F|^{-{{rB}\over 2}}\big)^{-s-{{A+2C}\over B}}
\int_{E\in {\cal M}_{{\bf A}_F,r},N(E)\leq |\Delta_F|^{r\over 2}}
N(E)^{\alpha\cdot(-s-{{A+2C}\over B})+\beta}d\mu(E).\cr}$$
 From the calculation in (2.2),  when
${\rm\,  Re}(\alpha\cdot s+\beta)<0$ and ${\rm\, 
Re}\big(\alpha\cdot(-s-{{A+2C}\over
B}\big)+\beta\big)>0$, being holomorphic functions,
$$II(s)=-W_F(r)\cdot {1\over {\alpha\cdot s+\beta}}\qquad{\rm\,  and}\qquad 
IV(s)=W_E(r)\cdot
{1\over {\alpha\cdot(-s-{{A+2C}\over B}\big)+\beta}}.$$
\vskip 0.30cm
\noindent
{\it Type 2}. Convergence for $I(s)$ and $I(-s-{{A+2C}\over B})$ where
$$I(s)=\big(|\Delta_F|^{-{{rB}\over 2}}\big)^s\int_{E\in {\cal
M}_{{\bf A}_F,r},N(E)\leq |\Delta_F|^{r\over 2}}\Big(H_{\rm\, 
ga}^0(F,E)^A\cdot
N(E)^{Bs+C}-N(E)^{\alpha s+\beta}\Big)d\mu(E).$$

By the discussion above about $II(s)$, unless $B=\alpha$, $I(s)$ and 
$I(-s-{{A+2C}\over B})$, and
hence $Z_{F,r;A,B,C;\alpha,\beta}(s)$ cannot be meromorphically extended as 
a meromorphic
function to the whole $s$-plane. (See also the discussion below.) Thus, we 
introduce the following
\vskip 0.30cm
\noindent
{\bf Compatibility Conditions}: $\alpha=B$ and $\beta=C$.
\vskip 0.30cm
With this,  by a change of variables, we also assume that
$B=1$ and $C=0$ so as to obtain the following integration:
$$Z_{F,r;A}(s):=Z_{F,r;A,-1,0;-1,0}(s):=\big(|\Delta_F|^{{{r}\over 2}}\big)^s
\cdot\int_{E\in {\cal M}_{{\bf A}_F,r}}\Big(H_{\rm\, 
ga}^0(F,E)^A-1\Big)\cdot N(E)^{-s}d\mu(E).$$

\noindent
Furthermore, for any $g\in {\rm\,  GL}_r({\bf A})$, from the definitionand 
the fact that
$H^0({\bf A}_F,g)$ is discrete, by writing down each term precisely,
$$H_{\rm\,  ga}^0(F,g)=1+\sum_{\alpha\in H^0({\rm\,  Spec}{\cal O}_F,{\cal 
E}(g))\backslash
\{0\}}\exp\Big(-\pi\sum_{\sigma:{\bf R}}|g_\sigma\cdot\alpha|^2-
2\pi\sum_{\sigma:{\bf C}}|g_\sigma\cdot \alpha|^2\Big)=:1+H_{\rm\, 
ga}'(F,g).$$ In this expression,
the first term is simply the constant function 1 on the moduli space, while 
each term in the second
decays exponentially. With this, by a standard argument about convergence 
ofan integration of theta
series for higher rank lattices in reduction theory, see. e.g., [25, Chap. 
III, Lect. 15], and the
fact that 1 in the first term $H_{\rm\,  ga}^0(F,E)^A$ cancels with 
thesecond term 1 in the
combination $H_{\rm\,  ga}^0(F,E)^A-1$, we conclude that
$II(s)$ and $II(-s-A)$ are all holomorphic functions, provided $A>0$. Allin 
all,
we have proved the following
\vskip 0.30cm
\noindent
{\bf Main Theorem.} {\it For any strictly
positive real number
$A$,
$$Z_{F,r;A}(s)=:\big(|\Delta_F|^{{{r}\over 2}}\big)^s\cdot
\int_{E\in {\cal M}_{{\bf A}_F,r}}\Big(H_{\rm\, \, ga}^0(F,E)^A-1\Big)\cdot 
N(E)^{-s}d\mu(E)$$
is holomorphic when ${\rm\, \, Re}(s)>A$. Moreover,

\noindent
(1) $Z_{F,r;A}(s)$ admits a meromorphic continuation to the whole complex 
$s$ plane which
only has simple poles at $s=0$ and $s=A$ with the same residue $W_F(r)$, 
i.e, the Tamagawa
volume of the moduli space ${\cal M}_{{\bf A}_F,r}[|\Delta_F|^{r\over 2}]$;

\noindent
(2) $Z_{F,r;A}(s)$ satisfies
the functional equation $$Z_{F,r;A}(s)=Z_{F,r;A}(A-s).$$}

\noindent
{\bf Main Definition.} {\it The function $Z_{F,r}(s):=Z_{F,r;1}(s)$ is 
called the rank $r$
non-abelian zeta function of $F$.}
\vskip 0.30cm
\noindent
{\it Remarks.} The latest definition may be justified by Iwasawa's ICM talk 
at MIT.
(See e.g., [13 and/or 31].) As they stand, our non-abelian zeta functions 
expose non-abelian
aspect of number fields. For details, see e.g., [31].

\noindent
(2) One may simply use moduli spaces ${\cal M}_{F,r}(d)$ to introduce new 
non-commutative zeta
functions for number fields as well.
\vskip 0.30cm
\noindent
{\bf Acknowledgement.} I would like to thank Ch. Deninger for the discussionand
introducing me the works of Stuhler and Grayson. Special thanks are also 
dueto G. van der Geer
for introducing and explaining to me his joint work with R. Schoof.
\vfill
\eject
\centerline {\li REFERENCES}
\vskip 0.40cm
\item{[1]} E. Bombieri \& J. Vaaler, On Siegel's lemma. Invent. Math.
{\bf 73} (1983), no. 1, 11--32.
\vskip 0.20cm
\item{[2]} A. Borel, Some finiteness properties of adele groups over
number fields, Publ. Math., IHES, {\bf 16} (1963) 5-30
\vskip 0.20cm
\item{[3]} A. Borisov, Convolution structures and arithmetic
cohomology, to appear in Comp. Math.
\vskip 0.20cm
\item{[4]} J.-B. Bost, Fibr\'es vectoriels hermitiens, degr\'e
d'Arakelov et polygones canoniques, Appendix A to  Exp. No. {\bf 795}, 
S\'eminaire
Bourbaki 1994/95, Ast\'erisque  {\bf 237} (1996), 154--161.
\vskip 0.25cm
\item{[5]} A. Connes, Trace formula in
noncommutative geometry and the zeros of the Riemann zeta function.  Sel.
math. New ser  {\bf 5}  (1999),  no. 1, 29--106.
\vskip 0.20cm
\item{[6]} Ch. Deninger,  On the $\Gamma$-factors attached to motives.
Invent. Math. {\bf 104} (1991), no. 2, 245--261.
\vskip 0.20cm
\item{[7]} Ch. Deninger, Motivic $L$-functions and regularized determinants,in
{\it Proc. Sympos. Pure Math}, {\bf 55}, {\it Motives}, edited by U. 
Jannsen, S. Kleiman and J.-P.
Serre,  (1994), 707-743
\vskip 0.25cm
\item{[8]} Ch. Deninger, Some analogies between number theory
and dynamical systems on foliated spaces. {\it Proceedings of the 
International Congress of
Mathematicians}, Vol. I (Berlin, 1998).  Doc. Math.  1998,  Extra Vol. I, 
163--186
\vskip 0.25cm
\item{[9]} G.B. Folland, {\it A course in abstract harmonic analysis}, 
Studies in advanced
mathematics, CRC Press, 1995
\vskip 0.20cm
\item{[10]} G. van der Geer \& R. Schoof, Effectivity of Arakelov
Divisors and the Theta Divisor of a Number Field, Sel. Math., New ser.
{\bf 6} (2000), 377-398
\vskip 0.20cm
\item{[11]} D. R. Grayson, Reduction theory using semistability.
Comment. Math. Helv.  {\bf 59}  (1984),  no. 4, 600--634.
\vskip 0.20cm
\item{[12]} D. R. Grayson,  Reduction theory
using semistability. II.  Comment. Math. Helv.  {\bf 61}  (1986),  no. 4,
661--676.
\vskip 0.20cm
\item{[13]} K. Iwasawa, Letter to Dieudonn\'e, April 8, 1952, in
  {\it Zeta Functions in Geometry}, edited by N.Kurokawa and T. Sunuda,
Advanced Studies in Pure Math. {\bf 21} (1992), 445-450
\vskip 0.20cm
\item{[14]} S. Lang, {\it Algebraic Number Theory},
Springer-Verlag, 1986
\vskip 0.20cm
\item{[15]} S. Lang, {\it Fundamentals on Diophantine Geometry},
Springer-Verlag, 1983
\vskip 0.20cm
\item{[16]} S. Lang, {\it Introduction to Arekelov Theory},
Springer-Verlag, 1988
\vskip 0.20cm
\item{[17]} X. Li,  A note on the Riemann-Roch theorem for
function fields.  {\it Progr. Math.}, 139, (1996),  567--570.
\vskip 0.25cm
\item{[18]} C. Moreno, {\it  Algebraic curves over finite fields.}
Cambridge Tracts in Mathematics, {\bf 97}, Cambridge University Press, 1991
\vskip 0.20cm
\item{[19]} A. Moriwaki, Stable sheaves on arithmetic curves, a personal
note dated in 1992
\vskip 0.20cm
\item{[20]} D. Mumford, {\it Geometric Invariant Theory}, Springer-Verlag,
(1965)
\vskip 0.20cm
\item{[21]} J. Neukirch, {\it Algebraic Number Theory}, Grundlehren der
Math. Wissenschaften, Vol. {\bf 322}, Springer-Verlag, 1999
\vskip 0.20cm
\item{[22]} A.N. Parshin, On the arithmetic of two-dimensional schemes. I. 
Distributions
and residues. (Russian)  Izv. Akad. Nauk SSSR Ser. Mat.  {\bf 40} (1976), 
no. 4, 736--773, 949.
\vskip 0.20cm
\item{[23]} J.-P. Serre,  Zeta and $L$ functions, in {\it Arithmetical 
Algebraic
Geometry} (Proc. Conf. Purdue Univ., 1963)   Harper \& Row (1965), 82-92
\vskip 0.25cm
\item{[24]} J.-P. Serre, {\it Algebraic Groups and Class Fields}, GTM
{\bf 117}, Springer-Verlag (1988)
\vskip 0.20cm
\item{[25]} C. L. Siegel, {\it Lectures on the geometry of numbers}, notes 
by B.
Friedman, rewritten by K. Chandrasekharan with the assistance of R. Suter, 
Springer-Verlag,
1989.
\vskip 0.20cm
\item{[26]} U. Stuhler, Eine Bemerkung zur Reduktionstheorie
quadratischer Formen.  Arch. Math.  {\bf 27} (1976), no. 6,
604--610.
\vskip 0.20cm
\item{[27]} U. Stuhler,  Zur Reduktionstheorie der positiven
quadratischen Formen. II.  Arch. Math.   {\bf 28}  (1977), no. 6, 611--619.
\vskip 0.20cm
\item{[28]} J. Tate, Fourier analysis in number fields and Hecke's
zeta functions, Thesis, Princeton University, 1950
\vskip 0.20cm
\item{[29]} A. Weil, {\it Basic Number Theory}, Springer-Verlag, 1973
\vskip 0.20cm
\item{[30]} L. Weng, Non-Abelian Class Field Theory for Riemann Surfaces, at
math.AG/0111240
\vskip 0.25cm
\item{[31]} L. Weng,  A Program for Geometric Arithmetic, at math.AG/0111241
\vskip 2.0cm
Lin WENG

Graduate School of Mathematics

Kyushu University

Fukuoka, 812-8581

JAPAN
\end